\documentclass{article}
\usepackage{amssymb}
\usepackage{amsmath}

\newtheorem{thm}{Theorem}[section]
\newtheorem{cor}[thm]{Corollary}
\newtheorem{lem}[thm]{Lemma}

\newtheorem{defn}[thm]{Definition}
\newtheorem{rem}[thm]{Remark}

\newcommand{\Real}{\mathbb R}

\newcommand{\half}{\frac{1}{2}}
\newcommand{\diffd}{\mathrm{d}} %integration and differentiation `d'
\newcommand{\Prob}{\mathbb{P}} %probability measure
\newcommand{\Proof}{\noindent\textbf{Proof: }}

\newcommand{\F}{\mathcal{F}} %for probability filtrations
\newcommand{\eL}{\mathcal{L}} %for L^1 function spaces etc.
\newcommand{\N}{\mathbb{N}} %the set of natural numbers
\newcommand{\Ind}{\mathbf{1}} %indicator function for sets

\newcommand{\cond}[2]{\bigl(#1\lvert#2\bigr)}
\newcommand{\ccond}[2]{\Bigl(#1\bigl\vert#2\Bigr)}
\newcommand{\seq}[1]{\bigl(#1\bigr)}  % round brackets
\newcommand{\bseq}[1]{\Bigl(#1\Bigr)}  % round brackets
\newcommand{\set}[1]{\bigl\{#1\bigr\}}
\newcommand{\bset}[1]{\Bigl\{#1\Bigr\}}
\newcommand{\bbset}[1]{\biggl\{#1\biggr\}}
\newcommand{\abs}[1]{\bigl\vert#1\bigr\vert}

\newcommand{\node}[1]{\text{node}_{#1}}
\newcommand{\Q}{\mathbb{Q}} %a new measure Q

\newcommand{\eg}{e.g.\ }
\newcommand{\qed}{\hfill$\Box$\bigskip}
\ifx\pdfoutput\undefined % We're not running pdftex
\else
\RequirePackage[colorlinks,hyperindex,plainpages=false]{hyperref}
\fi
\begin{document}
\title{A new formulation of the spine approach\\ to branching diffusions\thanks{
This \texttt{arXiv} article is an updated version of \emph{A new
formulation of the spine approach to branching diffusions}, (2004),
no.~0404, Mathematics Preprint, University of Bath.}}
\author{Robert Hardy and Simon C. Harris\thanks{Email: \texttt{S.C.Harris@bath.ac.uk}, Web: http://people.bath.ac.uk/massch}\\ University of Bath, UK }
\date{\today}
% MSC-class: 60J80
\maketitle

% ----------------------------------------------------------------
\begin{abstract}
We present a new formalization of the spine change of measure approach for branching diffusions
that improves on the scheme laid out for branching Brownian motion in
Kyprianou \cite{kyprianou:travelling_waves}, which itself made use of
earlier works of Lyons \emph{et al} \cite{lyons_pemantle_peres:conceptual_LlogL,
lyons:simple_biggins_convergence, lyons_kurtz_pemantle_peres:conceptual_Kesten_Stigum}.
We use our new formulation to interpret certain `Gibbs-Boltzmann' weightings of particles
and use this to give a new, intuitive and  proof of a more general
`Many-to-One' result which enables expectations of sums over particles in the branching diffusion to be
calculated purely in terms of an expectation of one particle.
Significantly, our formalization has provided the foundations that facilitate a variety of
new, greatly simplified and more intuitive proofs in branching diffusions: see, for example,
the $\eL^p$ convergence of additive martingales in Hardy and Harris \cite{L-p_convergence},
the path large deviation results for branching Brownian motion in Hardy and Harris \cite{bbm_large_deviations}
and the large deviations for a continuous-typed branching diffusion in Git \emph{et al} \cite{ghh}
and Hardy and Harris \cite{preprint-HW_large_deviations}.

\end{abstract}
% ----------------------------------------------------------------
\section{Introduction}\label{foundations_overview}

One of the central elements of the spine approach is to interpret
the behaviour of a branching process under a certain change of measure.
Such an interpretation was first laid out by Chauvin and Rouault
\cite{chauvin:KPP_app_to_spatial_trees} in the case of branching
Brownian motion, and we first briefly review the main ideas on a
heuristic level.

Consider a branching Brownian motion (BBM) with constant breeding
rate $r$, that is, a branching process whereby particles diffuse
independently according to a Brownian motion and at any moment
undergo fission at a rate $r$ to produce two particles that each
evolve independently from their birth position, and so on.
We suppose that the probabilities of this process are given by $\set{P^x: x \in \Real}$,
where $P^x$ is a measure defined on the natural filtration
$(\F_t)_{t \geq 0}$ such that it is the law of the process
initiated from a single particle positioned at $x$. Suppose that
the configuration of this branching Brownian motion at time $t$ is
given by the $\Real$-valued point process $\mathbb{X}_t :=
\set{X_u(t): u \in N_t}$ where $N_t$ is the set of individuals
alive at time $t$. It is well known that for any $\lambda \in
\Real$,
\begin{equation}\label{marky_mark}
Z_\lambda(t):=\sum_{u \in N_t} e^{-r t} e^{\lambda X_u(t)- \half
\lambda^2 t}
\end{equation}
defines a \emph{strictly-positive} $P$-martingale, so
$Z_\lambda(\infty) := \lim_{t \to \infty} Z_\lambda(t)$ is almost
surely finite under $P^x$. The important contribution of Chauvin
and Rouault \cite{chauvin:KPP_app_to_spatial_trees} was to
determine a pathwise construction of the measure $\Q_\lambda^x$
where
\begin{equation}\label{ken_dodd}
\frac{\diffd \Q_\lambda^x}{\diffd P^x}\bigg\vert_{\F_t} =
\frac{Z_\lambda(t)}{Z_\lambda(0)},
\end{equation}
with the term $Z_\lambda(0)$ acts as a normalizing factor. \bigbreak
The measure $\Q_\lambda^x$ defined at \eqref{ken_dodd} is equivalent
to the following pathwise construction:
\begin{itemize}
\item starting from position $x$, the original ancestor diffuses
according to a Brownian motion on $\Real$ with drift $\lambda$;

\item at an accelerated rate $2r$ the particle undergoes fission producing two
particles;

\item with equal probability, one of these two particles is
selected;

\item this chosen particle repeats stochastically the behaviour of
the parent;

\item the other particle initiates, from its birth position, an
independent copy of a $P^\cdot$ branching Brownian motion with
branching rate $r$.
\end{itemize}
%\end{defn}

The chosen line of descent in such pathwise constructions of the
measure, here $\Q_\lambda$, has come to be known as the \emph{spine}
as it can be thought of as the backbone of the branching process
$\mathbb{X}_t$ from which all particles are born. Although Chauvin
and Rouault's work on the measure change continued in a paper
co-authored with Wakolbinger \cite{chauvin:growing}, where the new
measure is interpreted as the result of building a conditioned tree
using the concepts of Palm measures, it wasn't until the so-called
`conceptual proofs' of Lyons, Kurtz, Peres and Pemantle published
around 1995 (\cite{lyons_pemantle_peres:conceptual_LlogL,
lyons:simple_biggins_convergence,
lyons_kurtz_pemantle_peres:conceptual_Kesten_Stigum}) that the spine
approach really began to crystalize. These papers laid out a formal
basis for spines using a series of new measures on two underlying
spaces of sample trees with and without distinguished lines of
descent (the spine). Of particular interest is the paper by Lyons
\cite{lyons:simple_biggins_convergence} which gave a spine-based
proof of the $\eL^1$-convergence of the well-known martingale for
the Galton-Watson process. Here we first saw the \emph{spine
decomposition} of the martingale as the key to using the intuition
provided by Chauvin and Rouault's pathwise construction of the new
measure -- Lyons used this together with a previously known
measure-theoretic result on Radon-Nikodym derivatives that allows us
to deduce the behaviour of the change-of-measure martingale under
the original measure by investigating its behaviour under the second
measure. Similar ideas have recently been used by Kyprianou
\cite{kyprianou:travelling_waves} to investigate the
$\eL^1$-convergence of the BBM martingale \eqref{marky_mark}, by
Biggins and Kyprianou \cite{bigginkypMTBP} for multi-type branching
processes in discrete time, by Geiger
\cite{geiger:elementary_GW_proofs, Geiger:PoissonGWtrees} for
Galton-Watson processes, by Georgii and Baake \cite{georgii} to
study ancestral type behaviour in a continuous time branching Markov
chain, as well as Olofsson \cite{olofsson:xlogx}. Also
%Lyons \emph{et al} \cite{lyons_kurtz_pemantle_peres:conceptual_Kesten_Stigum}
see Athreya \cite{athreya:bernoulli_change_of_measure_MC},
Geiger \cite{geiger:splitting_trees, Geiger:Sizebiasedsplittingtrees, Geiger:infvartree},
Iksanov \cite{iksanov}, Rouault and Liu \cite{rouault_liu:two_measures_on_tree_boundary}
and Waymire and Williams \cite{Waymire},
to name just a few other papers where
spine and size-biasing techniques have already proved extremely useful in branching process
situations.

\bigskip\noindent
In this article we present a new formalization of the spine
approach that improves the schemes originally laid out by Lyons
\emph{et al} \cite{lyons_pemantle_peres:conceptual_LlogL,
lyons:simple_biggins_convergence,
lyons_kurtz_pemantle_peres:conceptual_Kesten_Stigum} and later for BBM
in Kyprianou \cite{kyprianou:travelling_waves}.
Although the set-up costs of our spine formalization are quite large, at least in terms
of definitions and notation, the underlying ideas are all extremely simple and intuitive.
The power of this approach should not be underestimated.
The further techniques  developed subsequently in Hardy and Harris
\cite{L-p_convergence, bbm_large_deviations,preprint-HW_large_deviations},
Git \emph{et. al.} \cite{ghh} and J.W.Harris \& S.C.Harris \cite{hhx2} manage to
completely bypass many previous technical problems and difficult non-linear calculations,
with spine calculations facilitating the reduction
to relatively straightforward classical one-particle situations;
this paper serves as a foundation for these other works.

In the first instance our improvements correct a perceived weakness in the
Lyons \emph{et al} scheme where one of the measures they defined
had a time-dependent mass and could therefore not be normalized to
be a probability measure and lacked a clear interpretation in terms of any
direct process construction; an immediate consequence of this
improvement is that here \emph{all} measure changes are carried
out by \emph{martingales} and we regain a clear intuitive construction.

Another difference is in our use of filtrations and
sub-filtrations, where Lyons \emph{et al} instead used
marginalizing. As we shall show, this brings substantial benefits
since it allows us to relate the spine and the branching diffusion
through the conditional-expectation operation, and in this way
gives us a proper methodology for building \emph{new} martingales
for the branching diffusion based on known martingales for the
spine. In other work we have used this powerful construction to
obtain martingales for large-deviations problems, including a neat
proof of a large deviations principle for branching Brownian
motion in \cite{bbm_large_deviations},
and to study a BBM with quadratic branching rate in \cite{hhx2}.

The conditional-expectation approach also leads directly to new
and simple proofs of some key results for branching diffusions.
The first of these concerns the relation that becomes clear
between the spine and the `Gibbs-Boltzmann' weightings for the
branching particles. Such weightings are well known in the theory
of branching process, and Harris \cite{harris:gibbs_boltzmann}
contains some analysis for a model of a typed branching diffusion
that is similar in spirit to the models we shall be considering
here. In our formulation these important weightings can be
interpreted as a conditional expectation of a spine event, and we
can use them to immediately obtain a new and very useful
interpretation of the additive operations previously seen only
within the context of the Kesten-Stigum theorem and related
problems.
Another application of our approach we give a
substantially easier proof of a far more general %version of a
\emph{Many-to-One} theorem
%than that
that is so often useful in branching processes applications;
for example, in Champneys \emph{et al} \cite{art:aap} or
Harris and Williams \cite{art:largedevs1} it was a key tool in their more classical
approaches to branching diffusions.

The layout of this paper is as follows. In Section \ref{sec:models},
we will introduce the branching models; describing a binary
branching multi-type BBM that we will frequently use as an example,
before describing a more general branching Markov process model with
random family sizes. In Section \ref{sec:spaces}, we introduce the
\emph{spine} of the branching process as a distinguished infinite
line of descent starting at the initial ancestor, we describe the
underlying space for the branching Markov process with spine and we
also introduce various fundamental filtrations. In Section
\ref{measure_section},  we define some fundamental probability
spaces, including a probability measure for the branching process
with a randomly chosen spine. In Section  \ref{sec:martingales},
various martingales are introduced and discussed. In particular, we
see how to use filtrations and conditional expectation to build
`additive' martingales for the branching process out of the product
of three simpler `one-particle' martingales that only depend on the
behaviour along the path of the spine; used as changes of measure,
one martingale will increase the fission rate along the path of the
spine, another will size-bias the offspring distribution along the
spine, whilst the other one will change the motion of the spine.
Section \ref{sec:three_more_measures} discusses changes of measure
with these martingales and gives very important and useful intuitive
constructions for the branching process with spine under both the
original measure $\tilde{P}$ and the changed measure $\tilde{\Q}$.
Another extremely useful tool in the spine approach is the
\emph{spine decomposition} that we prove in Section
\ref{section:conditioning_on_the_spine}; this gives an expression
for the expectation of the `additive' martingale under the new
measure $\tilde\Q$ conditional on knowing the behaviour all along
the path of the spine (including the spine's motion, the times of
fission along the spine and number of offspring at each of the
spine's fissions). Finally, in Section \ref{sec:spine_methods}, we
use the spine formulation to derive an interpretation for certain
Gibbs-Boltzmann weights of $\tilde\Q$, discussing links with
theorems of Kesten-Stigum and Watanabe, in addition to proving a
very general `Many-to-One' theorem that enables expectations of sums
over particles in the branching process to be calculated as
expectation only involving the spine.

\section{Branching Markov models}\label{sec:models}
Before we present the underlying constructions for spines, it will
be useful to give the reader an idea of the branching-diffusion
models that we have in mind for applications. We first discuss a
finite-type branching diffusion, and then present a more general
model that shall be used as the basis of our spine constructions
in the following sections.

\subsection{A finite-type branching diffusion}\label{sec:finite-type_model}
Suppose that for a fixed $n \in \N$ we are given two sets of
positive constants $a(1), \ldots, a(n)$ and $r(1), \ldots, r(n)$.
Consider a typed branching diffusion in which the type of each
particle moves as a finite, irreducible and time-reversible Markov
chain on the set $I := \set{1, \ldots, n}$ with Q-matrix $\theta
Q$ ($\theta$ is a strictly positive constant that could be
considered as the \emph{temperature} of the system) and invariant
measure $\pi=\{\pi(1),\ldots,\pi(n)\}$. The spatial movement of a
particle of type $y$ is a driftless Brownian motion with
instantaneous variance $a(y)$, so that if $\seq{X_u(t),Y_u(t)} \in
\Real \times I$ is the space-type location of individual $u$ at
time $t$ then we have
\[
\diffd X_u(t) = a(Y_u(t)) \, \diffd B_t
\]
for a Brownian motion $B_t$. Fission of a particle of type $y$
occurs at a rate $r(y)$ to produce two particles at the same
space-type location as the parent.

We define $J:=\mathbb{R}\times I$, and suppose that the
configuration of this whole branching diffusion at time $t$ is
given by the $J$-valued point process $\mathbb{X}_t =
\set{\seq{X_u(t), Y_u(t)}: u \in N_t}$ where $N_t$ is the set of
individuals alive at time $t$. Let the measures $\set{P^{x,y}:
(x,y) \in \Real^2}$ on the filtered space $(\Omega, \F_\infty,
(\F_t)_{t \geq 0})$ be such that under $P^{x,y}$ a single initial
ancestor starting at $(x,y)$ evolves as the
branching diffusion, $\mathbb{X}_t$, as described above and where
$(\F_t)_{t \geq 0}$ the natural filtration generated by
the point process $\mathbb{X}_t$.

In this branching diffusion, each particle moves
in a stochastically similar manner:
let a process $(\xi_t,\eta_t)$ on $J$ under a measure $\Prob$
behaves stochastically like a single particle in the branching-diffusion
$\mathbb{X}_t$ with no branching occurring.
Thus, $\eta_t$ is an irreducible, time-reversible Markov
chain on $I$ with Q-matrix $\theta Q$ and invariant measure
$\pi=\{\pi_1,\ldots,\pi_n\}$, whilst $\xi_t$ moves as a driftless Brownian
motion and diffusion coefficient $a(y)>0$ whenever
$\eta_t$ is in state $y$:
\[
\diffd \xi_t = a(\eta_t)^\half \,\diffd B_t,
\]
for a $\Prob$-Brownian motion $B_t$. We note that the formal
generator of this process $(\xi_t,\eta_t)$ is:
\begin{equation}\label{eqn:SPM_generator}
\mathcal{H}F(x,y)=\half a(y)\frac{\partial^2 F}{\partial x^2} +
\theta \sum_{j\in I} Q(y,j) F(x,j),\qquad(F:J\to\Real).
\end{equation}
We shall often refer to the typed diffusion $(\xi_t, \eta_t) \in
\Real \times I$ as the \emph{single particle model}, after the
work carried out by Harris and Williams \cite{art:largedevs1}.

In Hardy \cite{hardy:thesis} and Hardy and Harris
\cite{L-p_convergence} this finite-type branching
diffusion has been investigated, and we briefly mention that the
proofs have been based on the following two martingales, the first
based on the whole branching diffusion and the second based only
on the single-particle model:
\begin{gather}
Z_\lambda(t) := \sum_{u \in N_t} v_\lambda(Y_u(t)) e^{\lambda
X_u(t) - E_\lambda t},\\
\zeta_\lambda(t) := e^{\int_0^t R(\eta_s) \, \diffd s}
v_\lambda(\eta_s) e^{\lambda \xi_t - E_\lambda t},\label{sp_mg}
\end{gather}
where $v_\lambda$ and $E_\lambda$ satisfy:
\[
\bseq{\half\lambda^2A+\theta Q+R} v_\lambda = E_\lambda v_\lambda,
\]
which is to say that $v_\lambda$ is an eigenvector of the matrix
$\half\lambda^2A+\theta Q+R$, with eigenvalue $E_\lambda$. These
two martingales should be compared with the corresponding
martingales \eqref{marky_mark} and $e^{\lambda B_t - \half
\lambda^2 t}$ for BBM and a single Brownian motion respectively.

\subsection{A general branching Markov process}

The spine constructions in our formulation can be applied to a
much more general branching Markov model, and we shall base the
presentation on the following model, where particles move
independently in a general space $J$ as a stochastic copy of some
given Markov process $\Xi_t$, and at a location-dependent rate
undergo fission to produce a location-dependent random number of
offspring that each carry on this branching behaviour
independently.
\begin{defn}[A general branching Markov process]\label{defn_gen_bp}
We suppose that three initial elements are given to us:
\begin{itemize}
\item a Markov process $\Xi_t$ in a measurable space $(J,
\mathcal{B})$,

\item a measurable function $R: J \to [0, \infty)$,

\item for each $x \in J$ we are given a random variable $A(x)$
whose probability distribution on the natural numbers $\set{0, 1,
\ldots}$ is $P\seq{A(x) = k} = p_k(x)$, and whose mean is $m(x) :=
\sum_{k = 0}^\infty kp_k(x)<\infty$.
\end{itemize}
From these ingredients we can build a branching process in $J$
according to the following recipe:
\begin{itemize}
\item Each particle of the branching process will live, move and
die in this space $(J, \mathcal{B})$, and if an individual $u$ is
alive at time $t$ we refer to its location in $J$ as $X_u(t)$.
Therefore the time-$t$ configuration of the branching process is a
$J$-valued point process $\mathbb{X}_t := \set{X_u(t): u \in N_t}$
where $N_t$ denotes the collection of all particles alive at time
$t$.

\item For each individual $u$, the stochastic behaviour of its
motion in $J$ is an independent copy of the given process $\Xi_t$.

\item The function $R: J \to [0,\infty)$ determines the rate at
which each particle dies: given that $u$ is alive at time $t$, its
probability of dying in the interval $[t, t+\diffd t)$ is
$R(X_u(t)) \diffd t + o(\diffd t)$.

\item If a particle $u$ dies at location $x \in J$ it is replaced
by $1 + A_u$ particles all positioned at $x$, where $A_u$ is an
independent copy of the random variable $A(x)$. All particles,
once born, progress independently of each other.
\end{itemize}
We suppose that the probabilities of this branching process are
$\set{P^x:x \in J}$ where under $P^x$ one initial ancestor
starts out at $x$.
\end{defn}
We shall first give a formal construction of the underlying
probability space, made up of the sample trees of the branching
process $\mathbb{X}_t$ in which the spines are the distinguished
lines of descent. Once built, this space will be filtered in a
natural way by the underlying family relationships of each sample
tree, the diffusing branching particles and the diffusing spine,
and then in section \ref{measure_section} we shall explain how we
can define new probability measures $\tilde{P}^x$ that extend each
$P^x$ up to the finest filtration that contains all information
about the spine and the branching particles.
Much of the notation that we use
for the underlying space of trees, the filtrations and the
measures is closely related to that found in Kyprianou
\cite{kyprianou:travelling_waves}.

Although we do not strive to present our spine approach in the
greatest possible generality, our general model already covers many
important situations whilst still being able to clearly demonstrate
all the key spine ideas. In particular, in all our models, new
offspring always inherit the position of their parent, although the
same spine methods should also readily adapt to situations with
random dispersal of offspring.

For greater clarity, we often use the finite-type
branching diffusion of Section \ref{sec:finite-type_model}
to introduce the ideas before following up
with the general formulation. For example, in this
finite-type model we would take the process $\Xi_t$ to be the
single-particle process $(\xi_t,\eta_t)$
which lives in the space $J := \Real \times I$
%described in the above section \ref{sec:finite-type_model},
and has generator $\mathcal{H}$ given by \eqref{eqn:SPM_generator}.
The birth rate in this model at location $(x,y)\in J$ will be
independent of $x$ and given by the function $R(y)$ for all $y \in
I$ and, since only binary branching occurs in this case, we also
have $P(A(x,y) = 1) = 1$ for all $(x,y) \in J$.

\section{The underlying space for spines}\label{sec:spaces}

\subsection{Marked Galton-Watson trees with spines} The set of
Ulam-Harris labels is to be equated with the set $\Omega$ of
finite sequences of strictly-positive integers:
\[
\Omega := \set{\emptyset} \cup \bigcup_{n\in\N} (\N)^n,
\]
where we take $\N=\set{1,2,\ldots}$. For two words $u,v\in\Omega$,
$uv$ denotes the concatenated word ($u\emptyset = \emptyset u =
u)$, and therefore $\Omega$ contains elements like `$213$' (or
`$\emptyset 213$'), which we read as `the individual being the 3rd
child of the 1st child of the 2nd child of the initial ancestor
$\emptyset$'. For two labels $v, u \in \Omega$ the notation $v <
u$ means that $v$ is an \emph{ancestor} of $u$, and $\abs{u}$
denotes the length of $u$. The set of all ancestors of $u$ is
equally given by
\[
\set{v:v<u} = \set{v: \exists w \in \Omega \text{ such that
}vw=u}.
\]

Collections of labels, ie. subsets of $\Omega$, will therefore be
groups of individuals. In particular, a subset $\tau\subset\Omega$
will be called a \emph{Galton-Watson tree} if:
\begin{enumerate}
\item $\emptyset\in\tau$,

\item if $u,v\in\Omega$, then $uv\in\tau$ implies $u\in\tau$,

\item for all $u\in\tau$, there exists $A_u\in{0,1,2,\ldots}$ such
that $uj\in\tau$ if and only if $1\leq j \leq 1+A_u$, (where
$j\in\N$).
\end{enumerate}
That is just to say that a Galton-Watson tree:
\begin{enumerate}
\item has a single initial ancestor $\emptyset$,

\item contains all ancestors of any of its individuals $v$,

\item has the $1+A_u$ children of an individual $u$ labelled in a
consecutive way,
\end{enumerate}
and is therefore just what we imagine by the picture of a family
tree descending from a single ancestor. Note that the `$1\leq j
\leq 1+A_u$' condition in 3 means that each individual has
\emph{at least} one child, so that in our model we are insisting
that Galton-Watson trees \emph{never die out}.

The set of all Galton-Watson trees will be called $\mathbb{T}$.
Typically we use the name $\tau$ for a particular tree, and whenever
possible we will use the letters $u$ or $v$ or $w$ to refer to the
labels in $\tau$, which we may also refer to as \emph{nodes of
$\tau$} or \emph{individuals in $\tau$} or just as \emph{particles}.

\bigskip\noindent
Each individual should have a \emph{location} in $J$ at each
moment of its \emph{lifetime}. Since a Galton-Watson tree
$\tau\in\mathbb{T}$ in itself can express only the \emph{family}
structure of the individuals in our branching random walk, in
order to give them these extra features we suppose that each
individual $u\in\tau$ has a mark $(X_u,\sigma_u)$ associated with
it which we read as:
\begin{itemize}
\item $\sigma_u\in\Real^+$ is the \emph{lifetime} of $u$, which
determines the \emph{fission time} of particle $u$ as $S_u :=
\sum_{v \leq u} \sigma_v$ (with $S_\emptyset:= \sigma_\emptyset$).
The times $S_u$ may also be referred to as the \emph{death} times;

\item $X_u:[S_u-\sigma_u,S_u)\to J$ gives the \emph{location} of
$u$ at time $t\in[S_u - \sigma_u,S_u)$.
\end{itemize}
To avoid ambiguity, it is always necessary to decide whether a
particle is in existence or not at its death time.
\begin{rem}\label{rem:death_and_birth_time_convention}
Our convention throughout will be that a particle $u$ dies `just
before' its death time $S_u$ (which explains why we have defined
$X_u:[S_u - \sigma_u,S_u)\to \cdot$ for example). Thus at the time
$S_u$ the particle $u$ has \emph{disappeared}, replaced by its
$1+A_u$ children which are all alive and ready to go.
\end{rem}

We denote a single marked tree by $(\tau,X,\sigma)$ or $(\tau, M)$
for shorthand, and the set of all marked Galton-Watson trees by
$\mathcal{T}$:
\begin{itemize}
\item $\mathcal{T} := \bbset{(\tau,X,\sigma):
\tau\in\mathbb{T}\text{ and for each }u\in\tau,
\sigma_u\in\Real^+, X_u:[S_u - \sigma_u,S_u)\to J}$.

\item For each $(\tau, X, \sigma) \in \mathcal{T}$, the set of
particles that are alive at time $t$ is defined as $N_t :=
\set{u\in\tau: S_u - \sigma_u \leq t < S_u}$.
\end{itemize}
Where we want to highlight the fact that these values depend on
the underlying marked tree we write \eg $N_t((\tau,X,\sigma))$ or
$S_u((\tau,M))$.

Any particle $u\in\tau$ that comes into existence creates a
\emph{subtree} made up from the collection of particles (and all
their marks) that have $u$ as an ancestor -- and $u$ is the
original ancestor of this subtree.
\begin{itemize}
\item $(\tau,X,\sigma)^u_j$, or $(\tau,M)^u_j$ for shorthand, is
defined as the \emph{subtree} growing from individual $u$'s $j$th
child $uj$, where $1 \leq j \leq 1+A_u$.
\end{itemize}
This subtree is a marked tree itself, but when considered as a
part of the original tree we have to remember that it comes into
existence at the space-time location $(X_u(S_u - \sigma_u),S_u - \sigma_u)$ --
which is just the space-time location of the death of particle $u$
(and therefore the space-time location of the birth of its child
$uj$).

\bigskip\noindent
Before moving on there is a further useful extension of the notation:
for any particle $u$ we extend the definition of $X_u$ from
the time interval $[S_u - \sigma_u, S_u)$ to allow all earlier
times $t\in[0,S_u)$:
\begin{defn}\label{defn:notation_extension}
Each particle $u$ is alive in the time interval $[S_u -
\sigma_u,S_u)$, but we extend the concept of its path in $J$ to
all earlier times $t < S_u$:
\[
X_u(t) := \left\{\begin{array}{cl}
  X_u(t) & \text{if } S_u - \sigma_u \leq t < S_u \\
  X_v(t) & \text{if } v< u \text{ and } J_v \leq t < S_v \\
\end{array}\right.
\]
\end{defn}
Thus particle $u$ inherits the path of its unique line of
ancestors, and this simple extension will allow us to later write
expressions like $\exp\{\int_0^t f(s)\,\diffd X_u(s)\}$ whenever $u
\in N_t$, without worrying about the birth time of $u$.

\bigskip\noindent
For any given marked tree $(\tau, M) \in \mathcal{T}$ we can
identify distinguished lines of descent from the initial ancestor:
$\emptyset, u_1, u_2, u_3,\ldots \in \tau$, in which $u_3$ is a
child of $u_2$, which itself is a child of $u_1$ which is a child
of the original ancestor $\emptyset$. We'll call such a subset of
$\tau$ a \emph{spine}, and will refer to it as $\xi$:
\begin{itemize}
\item a \emph{spine} $\xi$ is a subset of nodes $\set{\emptyset, u_1,
u_2, u_3,\ldots}$ in the tree $\tau$ that make up a unique line of
descent. We use $\xi_t$ to refer to the unique node in $\xi$ that
that is alive at time $t$.
\end{itemize}
In a more formal definition, which can for example be found in the
paper by Rouault and Liu
\cite{rouault_liu:two_measures_on_tree_boundary}, a spine is
thought of as a point on $\partial\tau$ the boundary of the tree
-- in fact the boundary is \emph{defined} as the set of all
infinite lines of descent. This explains the notation $\xi \in
\partial\tau$ in the following definition: we augment the space
$\mathcal{T}$ of marked trees to become
\begin{itemize}
\item $\tilde{\mathcal{T}} := \bset{(\tau,M, \xi):
(\tau,M)\in\mathcal{T}\text{ and }\xi\in \partial\tau}$ is the set
of \emph{marked trees with distinguished spines}.
\end{itemize}

\bigskip\noindent
It is natural to speak of the \emph{position of the spine at time
$t$} which think of just as the position of the unique node that
is in the spine and alive at time $t$:
\begin{itemize}
\item we define the time-$t$ position of the spine as $\xi_t :=
X_u(t)$, where $u\in \xi \cap N_t$.
\end{itemize}
By using the notation $\xi_t$ to refer to both the node in the
tree and that node's spatial position we are introducing potential
ambiguity, but in practice the context will make clear which we
intend. However, in case of needing to emphasize, we shall give
the node a longer name:
\begin{itemize}
\item $\node{t}((\tau,M,\xi)) := u$ if $u \in \xi$ is the node in
the spine alive at time $t$,
\end{itemize}
which may also be written as $\node{t}(\xi)$.

Finally, it will later be important to know how many fission times
there have been in the spine, or what is the same, to know which
generation of the family tree the node $\xi_t$ is in (where the
original ancestor $\emptyset$ is considered to be the 0th
generation)
\begin{defn}
We define the counting function
\[
n_t = \abs{\node{t}(\xi)},
\]
which tells us which generation the spine node is in, or
equivalently how many fission times there have been on the spine.
For example, if $\xi_t = \seq{\emptyset, u_1, u_2}$ then both
$\emptyset$ and $u_1$ have died and so $n_t = 2$.
\end{defn}

\subsection{Filtrations}

The reader who is already familiar with the Lyons \emph{et al}
\cite{lyons_kurtz_pemantle_peres:conceptual_Kesten_Stigum,lyons:simple_biggins_convergence,lyons_pemantle_peres:conceptual_LlogL}
papers will recall that they used two separate underlying spaces
of marked trees \emph{with} and \emph{without} the spines, then
marginalized out the spine when wanting to deal only with the
branching particles as a whole.
Instead, we are going to use the single underlying space $\tilde{\mathcal{T}}$, but define
\emph{four} filtrations of it that will encapsulate different knowledge.

\subsubsection{Filtration $(\F_t)_{t\geq0}$}
We define a filtration of $\tilde{\mathcal{T}}$ made up of the
$\sigma$-algebras:
\[
\F_t := \sigma \bseq{ (u,X_u,\sigma_u): S_u \leq t \,\, ; \,\,
(u,X_u(s):s \in [S_u - \sigma_u, t]): t \in [S_u - \sigma_u,S_u)},
\]
which in words means that $\F_t$ is generated by all the
information regarding the branching particles that have
\emph{lived and died} before time $t$ (this is the condition
$S_u\leq t$), along with just the information up to time $t$ of
those particles \emph{still alive} at time $t$ (this is the $t \in
[S_u - \sigma_u ,S_u)$ condition). Each of these $\sigma$-algebras
will be a subset of the limit defined as
\[
\F_\infty := \sigma\bseq{\bigcup_{t\geq0} \F_t}.
\]

\subsubsection{Filtration $(\tilde{\F}_t)_{t\geq0}$}\label{subsec:filt_F_tilde}

In order to know about the spine, we make this filtration finer,
defining $\tilde{\F}_t$ by adding into $\F_t$ the knowledge of
which node is the spine at time $t$:
\[
\tilde{\F}_t := \sigma\seq{\F_t, \node{t}(\xi)},\qquad
\tilde{\F}_\infty:=\sigma\bseq{\bigcup_{t\geq0} \tilde{\F_t}}.
\]
Consequently this filtration knows \emph{everything} about the
branching process and \emph{everything} about the spine: it knows
which nodes make up the spine, when they were born, when they died
(ie. the fission times $S_u$), and their family sizes.

\subsubsection{Filtration $(\mathcal{G}_t)_{t\geq0}$}

We define a filtration of $\tilde{\mathcal{T}}$,
$\set{\mathcal{G}_t}_{t\geq0}$, where the $\sigma$-algebras
\[
\mathcal{G}_t := \sigma\seq{\xi_s: 0\leq s\leq t},\qquad
\mathcal{G}_\infty:=\sigma\bseq{\bigcup_{t\geq0} \mathcal{G}_t},
\]
are generated by \emph{only} the spatial motion of the spine in
the $J$. Note that the events $G\in\mathcal{G}_t$ do not know
which \emph{nodes} of the tree $\tau$ actually make up the spine.

\subsubsection{Filtration $(\tilde{\mathcal{G}}_t)_{t\geq0}$}

We augment $\mathcal{G}_t$ by adding in information on the nodes
that make up the spine (as we did from $\F_t$ to $\tilde{\F}_t$),
as well as the knowledge of when the fission times occurred on the
spine and how big the families were that were produced:
\[
\tilde{\mathcal{G}}_t := \sigma\seq{\mathcal{G}_t, (\node{s}(\xi):s\leq
t), (A_u: u<\xi_t)},\qquad
\tilde{\mathcal{G}}_\infty:=\sigma\bseq{\bigcup_{t\geq0}
\tilde{\mathcal{G}_t}}.
\]

\subsubsection{Summary}
In brief, the key filtrations we shall make key use of are:
\begin{itemize}
\item $\F_t$ knows everything that has happened to all the
branching particles up to the time $t$, \emph{but does not know
which one is the spine};

\item $\tilde{\F}_t$ knows everything that $\F_t$ knows and also
knows which line of descent is the spine (it is in fact the finest
filtration);

\item $\mathcal{G}_t$ knows only about the spine's motion in $J$
up to time $t$, but does not actually know which line of descent
in the family tree makes up the spine;

\item $\tilde{\mathcal{G}}_t$ knows about the spine's motion and
also knows which nodes it is composed of. Furthermore it knows
about the fission times of these nodes and how many children were
born at each time.
\end{itemize}

\noindent
We note the obvious relationships between these
filtrations of $\tilde{\mathcal{T}}$ that
$\F_t \subset \tilde{\F}_t$ and
$\mathcal{G}_t  \subset \tilde{\mathcal{G}}_t  \subset \tilde{\F}_t$.
Trivially, we also note that $\mathcal{G}_t \nsubseteq \F_t$,
since the filtration $\F_t$ does not know \emph{which} line
of descent makes up the spine.

\section{Probability measures}\label{measure_section}
Having now carefully defined the underlying space for our
probabilities, we remind ourselves of the probability measures:
\begin{defn}
For each $x \in J$, let $P^x$ be the measure on
$(\tilde{\mathcal{T}},\F_\infty)$ such that the filtered
probability space $(\tilde{\mathcal{T}},\F_\infty,(\F_t)_{t \geq
0}, P^\cdot)$ is the canonical model for $\mathbb{X}_t$, the
branching Markov process described in Definition
\ref{defn_gen_bp}.
\end{defn}
For details of how the measures $P^x$ are formally constructed on
the underlying space of trees, we refer the reader to the work of
Neveu \cite{neveu:arbres} and Chauvin
\cite{chauvin:prod_mgs_stopping_lines, chauvin:bellman_harris}.
Note, we could equally think of $P^x$ as a measure on $({\mathcal{T}}, \F_\infty)$,
but it is convenient to use the enlarged sample space $\tilde{\mathcal{T}}$ for all
our measure spaces, varying only the filtrations.

Our spine approach relies first on building a measure $\tilde{P}^x$
under which the spine is a single genealogical line of descent
chosen from the underlying tree. If we are given a sample tree
$(\tau,M)$ for the branching process, it is easy to verify that, if
at each fission we make a uniform choice amongst the offspring to
decide which line of descent continues the spine $\xi$, when $u \in
\tau$ we have
\begin{equation}\label{uniform_spine_choice_definition}
\text{Prob}(u \in \xi) = \prod_{v < u} \frac{1}{1+A_v}.
\end{equation}
(In the binary-branching case, for example,
Prob$(A_v= 1) = 1$ and then $\text{Prob}(u \in \xi) =2^{-n_t}$.)
This simple observation is the key to our method for extending the
measures, and for this we make use of the following representation
found in Lyons \cite{lyons:simple_biggins_convergence}.
\begin{thm}
If $f$ is a $\tilde{\F}_t$-measurable function then we can write:
\begin{equation}\label{lyons_f_tilde_representation}
f = \sum_{u\in N_t} f_u \Ind_{(\xi_t = u)}
\end{equation}
where $f_u$ is $\F_t$-measurable.
\end{thm}
As a simple example of this, in the case of the finite-typed branching
diffusion of Section \ref{sec:finite-type_model}, such a
representation would be:
\begin{equation}\label{finite-type_lyons_decomp}
e^{\int_0^t R(\eta_s) \, \diffd s}v_\lambda(\eta_t) \, e^{\lambda
\xi_t - E_\lambda t} = \sum_{u \in N_t} e^{\int_0^t R(Y_u(s)) \,
\diffd s} v_\lambda(Y_u(t)) \, e^{\lambda X_u(t) - E_\lambda t} \,
\Ind_{(\xi_t = u)}.
\end{equation}
\begin{defn}\label{defn_of_P_tilde}
Given the measure $P^x$ on $(\tilde{\mathcal{T}}, \F_\infty)$ we
extend it to the probability measure $\tilde{P}^x$ on
$(\tilde{\mathcal{T}},\tilde{\F}_\infty)$ by defining
\begin{equation}\label{abcd}
\int_{\tilde{\mathcal{T}}} f \,\,\diffd \tilde{P}^x :=
\int_{\tilde{\mathcal{T}}} \sum_{u \in N_t} f_u \prod_{v < u}
\frac{1}{1+A_v} \,\,\diffd P^x,
\end{equation}
for each $f\in m\tilde{\F}_t$ with representation like
\eqref{lyons_f_tilde_representation}.
\end{defn}
The previous approach to spines, exemplified in Lyons
\cite{lyons:simple_biggins_convergence}, used the idea of
\emph{fibres} to get a measure analogous to our $\tilde{P}$ that
could measure the spine. However, a perceived weakness in this
approach was that the corresponding measure had time-dependent total
mass and could not be normalized to become a probability measure
with an intuitive construction, unlike our $\tilde P$. Our new idea
of using the down-weighting term of
\eqref{uniform_spine_choice_definition} in the definition of
$\tilde{P}$ is crucial in ensuring that we get a very natural
\emph{probability} measure (look ahead to Lemma
\ref{model_under_p_tilde}), and leads to the very useful situation
in which \emph{all} measure changes in our formulation are carried
out by \emph{martingales}.
\begin{thm}\label{p-tilde_is_extension}
This measure $\tilde{P}^x$ really is an extension of $P^x$ in that
$P = \tilde{P}\vert_{\F_\infty}$.
\end{thm}
\Proof If $f \in m\F_t$ then the representation
\eqref{lyons_f_tilde_representation} is trivial and therefore by
definition
\[
\int_{\tilde{\mathcal{T}}} f \,\,\diffd \tilde{P} =
\int_{\tilde{\mathcal{T}}} f \times \bseq{\sum_{u \in N_t}
\prod_{v < u} \frac{1}{1+A_v}} \,\,\diffd P.
\]
However, it can be shown that $\sum_{u \in N_t} \prod_{v < u}
\frac{1}{1+A_v} = 1$ by retracing the sum back through the lines
of ancestors to the original ancestor $\emptyset$, factoring out
the product terms as each generation is passed. Thus
\[
\int_{\tilde{\mathcal{T}}} f \,\diffd \tilde{P} =
\int_{\tilde{\mathcal{T}}} f \,\diffd P.
\]\qed

\begin{defn}
The filtered probability space
$(\tilde{\mathcal{T}},\tilde{\F}_\infty,(\tilde{\F}_t)_{t \geq 0},
\tilde{P})$ together with $(\mathbb{X}_t, \xi_t)$ will be referred to
as the \textbf{canonical model with spines}.
\end{defn}

\bigskip\noindent
In the single-particle model of section
\ref{sec:finite-type_model} we assumed the existence of a separate
measure $\Prob$ and a process $(\xi_t, \eta_t)$ that behaved
stochastically like a `typical' particle in the typed branching
diffusion $\mathbb{X}_t$. In our formalization the \emph{spine} is
exactly the single-particle model:
\begin{defn}\label{dave}
We define the measure $\Prob$ on $\seq{\tilde{\mathcal{T}},
\mathcal{G}_\infty}$ as the projection of $\tilde{P}$:
\[
\Prob\vert_{\mathcal{G}_t} := \tilde{P}\vert_{\mathcal{G}_t}.
\]
\end{defn}
Under the measure $\Prob$ the spine process $\xi_t$ has exactly
the same law as $\Xi_t$.
\begin{defn}
The filtered probability space
$(\tilde{\mathcal{T}},\mathcal{G}_\infty,(\mathcal{G}_t)_{t \geq
0}, \Prob)$ together with the spine process $\xi_t$ will be
referred to as the \textbf{single-particle model}.
\end{defn}

\section{Martingales}\label{sec:martingales}
Starting with the single Markov process $\Xi_t$ that lives in $(J,
\mathcal{B})$ we have built $(\mathbb{X}_t, \xi_t)$, a branching
Markov process with spines, in which the spine $\xi_t$ behaves
stochastically like the given $\Xi_t$. In this section we are
going to show how \emph{any} given martingale for the spine
$\zeta(t)$ leads to a corresponding additive martingale for the
whole branching model.

We have actually seen an example of this already. For the
finite-type model of section \ref{sec:finite-type_model} we met
two martingales:
\begin{gather}
Z_\lambda(t) := \sum_{u\in N_t} v_\lambda\seq{Y_u(t)} e^{\lambda
X_u(t) - E_\lambda t},\label{hay}\\
\zeta_\lambda(t) := e^{\int_0^t R(\eta_s) \,\diffd s}
v_\lambda(\eta_t) e^{\lambda \xi_t - E_\lambda t}.\label{bee}
\end{gather}
Just from their very form it has always been clear that they are
closely related. What we shall later be demonstrating in full generality
in this section is that the key to their relationship comes
through generalising the following $\tilde{\F}_t$-measurable martingale
for the multi-type BBM model:
\begin{defn}
We define an $\tilde{\F}_t$-measurable martingale:
\begin{equation}
\tilde{\zeta}_\lambda (t) := \prod_{u < \xi_t} (1+A_v) \times
v_\lambda(\eta_t) e^{\lambda \xi_t - E_\lambda t}.\label{mtmastermg}
\end{equation}
\end{defn}
An important result that we show in this article, in a more
general form, is that $Z_\lambda(t)$ and $\zeta_\lambda(t)$ are
just conditional expectations of this new
martingale $\tilde{\zeta}_\lambda$:
\begin{itemize}
\item $Z_\lambda(t) =
\tilde{P}\cond{\,\tilde{\zeta}_\lambda(t)\,}{\,\F_t}$,

\item $\zeta_\lambda(t) =
\tilde{P}\cond{\,\tilde{\zeta}_\lambda(t)\,}{\,\mathcal{G}_t}$.
\end{itemize}
We emphasize that this relationship does not appear to have previously been
formalized, and that it is only \emph{possible} because of our new
approach to the definition of $\tilde{P}$ as a \emph{probability}
measure, and of our using filtrations to capture the different
knowledge generated by the spine and the branching particles.

Furthermore, in the general form that we present below it provides
a consistent methodology for using well-known martingales for a
single process $\xi_t$ to get new additive martingales for the
related branching process. In Hardy and Harris
\cite{bbm_large_deviations,preprint-HW_large_deviations}
we use these powerful ideas to
give substantially easier proofs of large-deviations problems in
branching diffusions than have previously been possible.
\bigbreak

\bigskip
Suppose that $\zeta(t)$ is a $\seq{\tilde{\mathcal{T}},
(\mathcal{G}_t)_{t \geq 0}, \tilde{P}}$-martingale, which is to
say that it is a $\mathcal{G}_t$-measurable function that is a
martingale with respect to the measure $\tilde{P}$. For example, in
the case of our finite-type branching diffusion this could be the
martingale $\zeta_\lambda(t)$ which is $\mathcal{G}_t$-measurable
since it refers only to the spine process $(\xi_t, \eta_t)$.
\begin{defn}
We shall call $\zeta(t)$ a \textbf{single-particle martingale},
since it is $\mathcal{G}_t$-measurable and thus depends only to the spine $\xi$.
\end{defn}
Any such single-particle martingale can be used to define an
additive martingale for the whole branching process via the
representation \eqref{lyons_f_tilde_representation}:
\begin{defn}\label{defn_of_additive_martingale}
Suppose that we can represent the martingale $\zeta(t)$ as
\begin{equation}\label{dee}
\zeta(t) = \sum_{u \in N_t} \zeta_u(t) \Ind_{(\xi_t = u)},
\end{equation}
for $\zeta_u(t) \in m\F_t$, as at
\eqref{lyons_f_tilde_representation}. We can then define an
$\F_t$-measurable process $Z(t)$ as
\[
Z(t) := \sum_{u \in N_t} e^{-\int_0^t m(X_u(s)) R(X_u(s)) \,
\diffd s} \zeta_u(t),
\]
and refer to $Z(t)$ as the \textbf{branching-particle martingale}.
\end{defn}
The martingale property $Z(t)$ will be established in Lemma \ref{z_and_zeta_from_zeta_tilde}
after first building another martingale, $\tilde\zeta(t)$, from the single-particle martingale $\zeta(t)$.
First, for clarity, we take a moment to discuss this definition of the
additive martingale and the terms like $\zeta_u(t)$.

If we return to our familiar martingales \eqref{hay} and
\eqref{bee}, it is clear that
\begin{equation}\label{see}
\zeta_\lambda(t)
= e^{\int_0^t R(\eta_s) \,\diffd s}
v_\lambda(\eta_t) e^{\lambda \xi_t - E_\lambda t}
= \sum_{u\in N_t} e^{\int_0^t R(Y_u(s)) \, \diffd
s} v_\lambda\seq{Y_u(t)} e^{\lambda X_u(t) - E_\lambda t} \,
\Ind_{(\xi_t = u)}.
\end{equation}
The `$\zeta_u$' terms of \eqref{dee} could be here replaced with a
more descriptive notation $\zeta_\lambda[(X_u,Y_u)](t)$, where
\[
\zeta_u(t) = \zeta_\lambda[(X_u, Y_u)](t) := e^{\int_0^t R(Y_u(s))
\, \diffd s} v_\lambda\seq{Y_u(t)} e^{\lambda X_u(t) - E_\lambda
t},
\]
can be seen to essentially be a functional of the space-type path
$(X_u(t), Y_u(t))$ of particle $u$. In this way the original
single-particle martingale $\zeta_\lambda$ would be understood as
a functional of the space-type path $(\xi_t, \eta_t)$ of the spine
itself and we could write
\[
\zeta_\lambda(t) = \zeta_\lambda[(\xi , \eta)](t) = \sum_{u \in
N_t} \zeta_\lambda[(X_u, Y_u)](t) \, \Ind_{(\xi_t = u)}.
\]
This is the idea behind the representation \eqref{dee}, and in
those typical cases where the single-particle martingale is
essentially a functional of the paths of the spine $\xi_t$, as is
the case for our $\zeta_\lambda(t)$, we should just think of
$\zeta_u$ as being that same functional but evaluated over the
path $X_u(t)$ of particle $u$ rather than the spine $\xi_t$. The
representation \eqref{dee} can also be used as a more general way of
treating other %single-particle
martingales that perhaps are not such a simple functional of the spine path.

Finally, from \eqref{see} it is clear that the additive martingale
being defined by definition \ref{defn_of_additive_martingale} is
our familiar $Z_\lambda(t)$:
\begin{align*}
Z_\lambda(t) &= \sum_{u \in N_t} e^{-\int_0^t R(Y_u(s)) \, \diffd s}
\zeta_\lambda[(X_u, Y_u)](t)
= \sum_{u\in N_t} v_\lambda\seq{Y_u(t)} e^{\lambda
X_u(t) - E_\lambda t}.
\end{align*}

Although definition \ref{defn_of_additive_martingale} will work in
general, in the main the spine approach is interested in
martingales that can act as Radon-Nikodym derivatives between
probability measures, and therefore we suppose from now on that
$\zeta(t)$ is \emph{strictly positive}, and therefore that the
additive martingale $Z(t)$ is strictly positive.
\bigbreak

\bigskip\noindent
The work of Lyons \emph{et al}
\cite{lyons:simple_biggins_convergence,
lyons_kurtz_pemantle_peres:conceptual_Kesten_Stigum,
lyons_pemantle_peres:conceptual_LlogL}, that of Chauvin and
Rouault \cite{chauvin:KPP_app_to_spatial_trees} and more recently
of Kyprianou \cite{kyprianou:travelling_waves} suggests that when
a change of measure is carried out with a branching-diffusion
additive martingale like $Z(t)$ it is typical to expect three
changes: the spine will gain a drift, its fission times will be
increased and the distribution of its family sizes will be
size-biased. In section \ref{sec:understanding_Q} we shall confirm
this, but we first take a separate look at the martingales that could
perform these changes, and which we shall combine to obtain a
martingale $\tilde{\zeta}(t)$ that will ultimately be used to
change the measure $\tilde{P}$.
\begin{thm}
The expression
\[
\prod_{v < \xi_t} \seq{1+ m(\xi_{S_v})} e^{-\int_0^t
m(\xi_s)R(\xi_s) \, \diffd s}
\]
is a $\tilde{P}$-martingale that will increase the rate at which
fission times occur along the spine from $R(\xi_t)$ to
$(1+m(\xi_t))R(\xi_t)$:
\[
\frac{\diffd \mathbb{L}^{((1+m(\xi))R(\xi)) }_t}
{\diffd \mathbb{L}^{(R(\xi))}_t} = \prod_{v < \xi_t} \seq{1+
m(\xi_{S_v})} e^{-\int_0^t m(\xi_s)R(\xi_s) \, \diffd s}
\]
where $\mathbb{L}^{(R(\xi))}$ is the law of the Poisson (Cox)
process with rate $R(\xi_t)$ at time $t$.
\end{thm}
\begin{thm}
The term
\[
\prod_{v < \xi_t} \frac{1+A_v}{1+ m(\xi_{S_v})}
\]
is a $\tilde{P}$-martingale that will change the measure by
size-biasing the family sizes born from the spine:
\[
\text{if }v< \xi_t, \text{ then}\qquad \text{Prob}(A_v = k) =
\frac{(1+ k)p_k(\xi_{S_v})}{1 + m(\xi_{S_v})}.
\]
\end{thm}
The product of these two martingales with the single-particle
martingale $\zeta(t)$ will simultaneously perform the three
changes mentioned above:
\begin{defn}\label{sp_mg_tilde}
We define a $\tilde{\F}_t$-measurable martingale as
\begin{align}
\tilde{\zeta}(t) &:= \prod_{v < \xi_t} (1+A_v)
e^{-\int_0^t m(\xi_s)R(\xi_s) \, \diffd s} \times \zeta(t)\notag\\
&= \prod_{v < \xi_t} \frac{1+A_v}{1+ m(\xi_{S_v})} \times \prod_{v
< \xi_t} \seq{1+ m(\xi_{S_v})} e^{-m\int_0^t R(\xi_s) \, \diffd s}
\times \zeta(t).\label{eqn:sp_mg_tilde}
% &=
%\prod_{u < \xi_t} (1+A_v) v_\lambda(\eta_t) e^{\lambda \xi_t -
%E_\lambda t}.
\end{align}
\end{defn}
Significantly, \emph{only} the motion of the spine and the behaviour along the immediate
path of the spine will be affected by any change of measure using this martingale.
Also note, this martingale is the general form of $\tilde{\zeta}_\lambda(t)$
that we defined at \eqref{mtmastermg} for our finite-type model.

The real importance of the size-biasing and fission-time-increase
operations is that they introduce the correct terms into
$\tilde{\zeta}(t)$ so that the following key relationships hold:
\begin{lem}\label{z_and_zeta_from_zeta_tilde}
Both $Z(t)$ and $\zeta(t)$ are projections of $\tilde{\zeta}(t)$
onto their filtrations: for all $t$,
\begin{itemize}
\item $Z(t) = \tilde{P}\cond{\,\tilde{\zeta}(t)\,}{\,\F_t}$,

\item $\zeta(t) =
\tilde{P}\cond{\,\tilde{\zeta}(t)\,}{\,\mathcal{G}_t}$.
\end{itemize}
\end{lem}
\bigbreak
\Proof We use the representation
\eqref{lyons_f_tilde_representation} of $\tilde{\zeta}(t)$:
\begin{equation}
\tilde{\zeta}(t) = \sum_{u \in N_t}\prod_{v < u} (1+A_v)
e^{-\int_0^t m(X_u(s))R(X_u(s)) \, \diffd s}  \zeta_u(t)
\Ind_{(\xi_t=u)}.
\end{equation}
From this it follows that
\begin{align*}
\tilde{P}\cond{\tilde{\zeta}(t)}{\F_t} &= \sum_{u \in N_t}
e^{-\int_0^t m(X_u(s))R(X_u(s)) \, \diffd s}
\zeta_u(t) \times \prod_{v < u} (1+A_v) \, \tilde{P}\cond{\Ind_{(\xi_t=u)}}{\F_t} \\
&= \sum_{u \in N_t} e^{-\int_0^t m(X_u(s))R(X_u(s)) \, \diffd s}
\zeta_u(t)  = Z(t),
\end{align*}
since $\tilde{P}\cond{\Ind_{(\xi_t=u)}}{\F_t} = \Ind_{(u\in N_t)}
\times \prod_{v < u} (1+A_v)^{-1}$.

On the other hand, the martingale terms in \eqref{eqn:sp_mg_tilde}
imply
\begin{align*}
\tilde{P}\cond{\tilde{\zeta}(t)}{\mathcal{G}_t} &= \zeta(t) \times
\tilde{P} \ccond{\prod_{v < \xi_t} (1+A_v) e^{-\int_0^t
m(\xi_s)R(\xi_s) \, \diffd s}}{\mathcal{G}_t} = \zeta(t),
\end{align*}
\qed

\section{Changing the measures}\label{sec:three_more_measures}
For the finite type model,  the single-particle martingale
$\zeta_\lambda(t)$ defined at \eqref{sp_mg} can be used to define a new
measure for the single-particle model (as in \cite{hardy:thesis}), via
\[
\frac{\diffd \Prob_\lambda}{\diffd
\Prob}\bigg\vert_{\mathcal{G}_t} =
\frac{\zeta_\lambda(t)}{\zeta(0)}.
\]
We have now seen the close relationships between the three martingales
$\zeta_\lambda$, $Z_\lambda$ and $\tilde{\zeta}_\lambda$:
\[
Z_\lambda(t) =
\tilde{P}\cond{\,\tilde{\zeta}_\lambda(t)\,}{\,\F_t}, \qquad
\zeta_\lambda(t) =
\tilde{P}\cond{\,\tilde{\zeta}_\lambda(t)\,}{\,\mathcal{G}_t},
\]
and in this section we show in a more general form how these close
relationships mean that a new measure $\tilde{\Q}_\lambda$ defined
in terms of $\tilde{P}$ as
\[
\frac{\diffd \tilde{\Q}_\lambda}{\diffd
\tilde{P}}\bigg\vert_{\tilde{\F}_t} =
\frac{\tilde{\zeta}_\lambda(t)}{\tilde{\zeta}_\lambda(0)},
\]
will induce measure changes on the sub-filtrations $\mathcal{G}_t$
and $\F_t$ of $\tilde{\F}_t$ whose Radon-Nikodym derivatives are
given by $\zeta_\lambda(t)$ and $Z_\lambda(t)$ respectively.
We will also give an extremely useful and intuitive construction of the
measures $\tilde{P}$ and $\tilde\Q$.

\bigskip\noindent
We recall that in our set up we have a finest filtration
$(\tilde{\F}_t)_{t \geq 0}$ associated with the measure
$\tilde{P}$, and two sub-filtrations $(\F_t)_{t \geq 0}$ with
measure $P$ and $(\mathcal{G}_t)_{t \geq 0}$ with measure $\Prob$.
The martingale $\tilde{\zeta}$ can change the measure $\tilde{P}$:
\begin{defn}
A measure $\tilde{\Q}$ on
$(\tilde{\mathcal{T}},\tilde{\F}_\infty)$ is defined via its
Radon-Nikodym derivative with respect to $\tilde{P}$:
\[
\frac{\diffd \tilde{\Q}}{\diffd
\tilde{P}}\bigg\vert_{\tilde{\F}_t} =
\frac{\tilde{\zeta}(t)}{\tilde{\zeta}(0)}.
\]
Recall, this notation means that for each event $F \in
\tilde{\F}_t$ we define $\tilde{\Q}(F) := \tilde{P}\seq{\Ind_F
Z(t)}$.
\end{defn}
As we did for the measures $P$ and $\Prob$ in Section \ref{measure_section}, we can restrict
$\tilde{\Q}$ to the sub-filtrations:
\begin{defn}
We define the measure $\Q$ on
$(\tilde{\mathcal{T}},\F_\infty,(\F_t)_{t \geq 0})$ via
\[
\Q := \tilde{\Q}\vert_{\F_\infty}.
\]
\end{defn}
\begin{defn}\label{dave2}
We define the measure $\hat{\Prob}$ on
$(\tilde{\mathcal{T}},\mathcal{G}_\infty,({\mathcal{G}}_t)_{t
\geq 0})$ via
\[
\hat{\Prob} := \tilde{\Q}\vert_{\mathcal{G}_\infty}.
\]
\end{defn}

A consequence of our new formulation in terms of filtrations and
the equalities of Lemma \ref{z_and_zeta_from_zeta_tilde} is that
the changes of measure are carried out by $Z(t)$ and $\zeta(t)$ on
their subfiltrations:
\begin{thm}\label{mark_riley}
\[
\frac{\diffd \Q}{\diffd P}\bigg\vert_{\F_t} = \frac{Z(t)}{Z(0)},
\qquad \text{and} \qquad \frac{\diffd \hat{\Prob}}{\diffd
\Prob}\bigg\vert_{\mathcal{G}_t} = \frac{\zeta(t)}{\zeta(0)}.
\]
\end{thm}
\Proof These two results actually follow from a more general
observation that if $\tilde{\mu}_1$ and $\tilde{\mu}_2$ are two
measures defined on a measure space $(\Omega,
\tilde{\mathcal{S}})$ with Radon-Nikodym derivative
\begin{equation}\label{liable}
\frac{\diffd \tilde{\mu}_2}{\diffd \tilde{\mu}_1} = f,
\end{equation}
and if $\mathcal{S}$ is a sub-$\sigma$-algebra of
$\tilde{\mathcal{S}}$, then the two measures $\mu_1:=
\tilde{\mu}_1\vert_{\mathcal{S}}$ and $\mu_2:=
\tilde{\mu}_2\vert_{\mathcal{S}}$ on $(\Omega,\mathcal{S})$ are
related by the conditional expectation operation:
\[
\frac{\diffd \mu_2}{\diffd \mu_1} = \tilde{\mu}_1
\cond{f}{\mathcal{S}}.
\]
The proof of this is that if $g\in m\mathcal{S}$ and
$S\in\mathcal{S}$ then
\begin{align*}
\int_S g \,\diffd\mu_2 &= \int_S g \, \diffd \tilde{\mu}_2
&\text{since $g$ is also in
$m\tilde{\mathcal{S}}$, and $S\in\mathcal{S}$ too,} \\
&= \int_S g \, f\,\diffd \tilde{\mu}_1 &\text{by \eqref{liable},}\\
&= \int_S \tilde{\mu}_1\cond{g f}{\mathcal{S}} \,\diffd \tilde{\mu}_1 &\text{by definition of the conditional expectation,}\\
&= \int_S g \,\tilde{\mu}_1\cond{f}{\mathcal{S}} \,\diffd
\tilde{\mu}_1 &\text{since $g$ is $\mathcal{S}$-measurable,}\\
& = \int_S g \,\tilde{\mu}_1\cond{f}{\mathcal{S}} \,\diffd \mu_1
 &\text{since everything is in $m\mathcal{S}$.}
\end{align*}
Applying this general result \eqref{liable} using the
relationships between the general martingales given in Lemma
\ref{z_and_zeta_from_zeta_tilde} concludes the proof. \qed

\subsection{Understanding the measure $\tilde{\Q}$}\label{sec:understanding_Q}

As the name suggests, we should be able to think of the spine as the
backbone of the branching process. This is made precise by the
following decomposition:
\begin{thm}
The measure $\tilde{P}$ on $\tilde{\F}_t$ can be decomposed as:
\begin{equation}\label{decomp_of_p_tilde}
\diffd \tilde{P} (\tau, M, \xi) = \diffd \Prob(\xi) \diffd
\mathbb{L}^{(R(\xi))}(n) \prod_{v < \xi_t} \frac{1}{1+A_v}
\prod_{v < \xi_t} p_{A_v}(\xi_{S_v}) \prod_{j=1}^{A_v} \diffd
P\seq{(\tau,M)^v_j},
\end{equation}
where $\mathbb{L}^{(R(\xi))}$ is the law of the Poisson (Cox)
process with rate $R(\xi_t)$ at time $t$, and we remember that
$n_t$ counts the number of fission times on the spine before time
$t$.
\end{thm}
We can offer a clear intuitive picture of this decomposition, which we summarize
in the following lemma.
\begin{lem}\label{model_under_p_tilde}
The decomposition of measure $\tilde{P}$ at \eqref{decomp_of_p_tilde}
enables the following construction:
\begin{enumerate}
\item the spine's motion is determined by the single-particle
measure $\Prob$;

\item the spine undergoes fission at time $t$ at rate $R(\xi_t)$;

\item at the fission time of node $v$ on the spine, the single spine particle is
replaced by $1+A_v$ children, with $A_v$ being chosen
independently and distributed according to the location-dependent
random variable $A(\xi_{S_v})$ with probabilities $(p_k(\xi_{S_v}):k=0,1,\ldots)$;

\item the spine is chosen uniformly from the $1+A_v$ children at the fission point $v$;

\item each of the remaining $A_v$ children gives rise to the independent subtrees $(\tau,
M)^v_j$, for $1 \leq j \leq A_v$, which are not part of the spine
and which are each determined by an independent copy of the
original measure $P$ shifted to their point and time of creation.
\end{enumerate}
\end{lem}

\bigskip\noindent
This decomposition of $\tilde{P}_t$ given at
\eqref{decomp_of_p_tilde} will allow us to interpret the measure
$\tilde{\Q}$ if we appropriately factor the components of the
change-of-measure martingale $\tilde{\zeta}(t)$ across this
representation. On $\tilde{\F}_t$,
\begin{align}
\diffd \tilde{\Q} &= \tilde{\zeta}(t) \,\diffd
\tilde{P}\notag\\
&= \zeta(t) \times e^{-\int_0^t R(\xi_s)\diffd s} \seq{1 +
m(\xi_{S_v})}^{n_t} \times \prod_{v < \xi_t} \frac{1+A_v}{1 +
m(\xi_{S_v})} \times
\diffd \tilde{P}\notag\\
&= \diffd \hat{\Prob}(\xi) \,\diffd
\mathbb{L}^{((1+m(\xi))R(\xi))}(n) \,\prod_{v < \xi_t}
\frac{1}{1+A_v} \prod_{v < \xi_t}\frac{1+A_v}{1+
m(\xi_{S_v})}
%p_{A_v}
p_{A_v}(\xi_{S_v})
\prod_{j=1}^{A_v} \diffd
P\seq{(\tau,M)^v_j}.\label{q_decomp}
\end{align}
Just as we did for $\tilde{P}$, we can offer a clear interpretation of this
decomposition:
\begin{lem}\label{model_under_q}
Under the measure $\tilde{\Q}$,
\begin{enumerate}
\item the spine process $\xi_t$ moves as if under the
\emph{changed} measure $\hat{\Prob}$;

\item the fission times along the spine occur
at an \emph{accelerated} rate $(1 + m(\xi_t))R(\xi_t)$;

\item at the fission time of node $v$ on the spine,
the single spine particle is replaced by $1+{A}_v$
children, with ${A}_v$ being chosen as an independent copy
of the random variable $\tilde{A}(y)$
which has the \emph{size biased} offspring distribution $((1+k)p_k(y)/(1+m(y)):k=0,1,\ldots)$,
where $y=\xi_{S_v}\in J$ is the spine's location at the time of fission;

\item the spine is chosen uniformly from the $1+{A}_v$ particles at the fission point $v$;

\item each of the remaining ${A}_v$ children gives rise to the independent subtrees
$(\tau,M)^v_j$, for $1 \leq j \leq {A}_v$, which are not part of the spine
and evolve as independent processes determined by the measure $P$ shifted to
their point and time of creation.
\end{enumerate}
\end{lem}
Such an interpretation of the measure $\tilde{\Q}$ was first given
by Chauvin and Rouault \cite{chauvin:KPP_app_to_spatial_trees} in
the context of BBM, allowing them to come to the important
conclusion that under the new measure $\Q$ the branching diffusion
remains largely unaffected, except that the Brownian particles of
a single (random) line of descent in the family tree are given a
changed motion, with an accelerated birth rate -- although they did not
have random family sizes, so the size-biasing aspect was not seen.
In the context of spines, size-biasing was first introduced in the
Lyons \emph{et al} papers \cite{lyons:simple_biggins_convergence,
lyons_kurtz_pemantle_peres:conceptual_Kesten_Stigum,
lyons_pemantle_peres:conceptual_LlogL}.
Kyprianou \cite{kyprianou:travelling_waves} presented the decomposition of equation
\eqref{q_decomp} and the construction of $\Q$ at Lemma \ref{model_under_q} for BBM with random family sizes,
but did not follow our natural approach starting with the \emph{probability} measure $\tilde{P}$
that has subtly facilitated various benefits.

\section{The spine
decomposition}\label{section:conditioning_on_the_spine}

One of the most important results introduced in Lyons
\cite{lyons:simple_biggins_convergence} was the so-called
\emph{spine decomposition}, which in the case of the additive martingale
\[
Z_\lambda(t) = \sum_{u \in N_t} v_\lambda(Y_u(t)) e^{\lambda
X_u(t) - E_\lambda t},
\]
from the finite-type branching diffusion would be:
\begin{equation}\label{hand}
\tilde{\Q}_\lambda \cond{Z_\lambda(t)}{\tilde{\mathcal{G}}_\infty}
= \sum_{u < N_t} v_\lambda(\eta_{S_u}) e^{\lambda \xi_{S_u} -
E_\lambda S_u} + v_\lambda(\eta_t) e^{\lambda \xi_t - E_\lambda
t}.
\end{equation}
To prove this we start by decomposing the martingale as
\[
Z_\lambda(t) = \sum_{u \in N_t, u \notin \xi} v_\lambda(Y_u(t))
e^{\lambda X_u(t) - E_\lambda t} + v_\lambda(\eta_t) e^{\lambda
\xi_t - E_\lambda t},
\]
which is clearly true since one of the particles $u \in N_t$ must
be the in the line of descent that makes up the spine $\xi$.
Recalling that the $\sigma$-algebra $\tilde{\mathcal{G}}_\infty$
contains all information about the line of nodes that makes up the
spine, all about the spine diffusion $(\xi_t, \eta_t)$ for all
times $t$, and also contains all information regarding the fission
times and number of offspring along the spine, it is useful to partition the particles $v \in
\set{u \in N_t, u \notin \xi}$ into the distinct subtrees
$(\tau,M)^u$ that were born at the fission times $S_u$ from the
particles that made up the spine before time $t$, or in other
words those nodes in the $\set{u < \xi_t}$ of ancestors of the
current spine node $\xi_t$. Thus:
\[
Z_\lambda(t) = \sum_{u < \xi_t} e^{\lambda \xi_{S_u} - E_\lambda
S_u} \bbset{\sum_{v \in N_t, v \in (\tau,M)^u} v_\lambda(Y_v(t))
e^{\lambda (X_u(t) - \xi_{S_u}) - E_\lambda (t-S_u)}} +
v_\lambda(\eta_t) e^{\lambda \xi_t - E_\lambda t}.
\]
If we now take the $\tilde{\Q}_\lambda$-conditional expectation of
this, we find
\begin{multline*}
\tilde{\Q}_\lambda\cond{Z_\lambda(t)}{\tilde{\mathcal{G}}_\infty}
= v_\lambda(\eta_t) e^{\lambda \xi_t-E_\lambda t} + \\
\sum_{u < \xi_t} e^{\lambda \xi_{S_u}-E_\lambda S_u}
\tilde{\Q}_\lambda\ccond{\sum_{v \in N_t, v \in (\tau,M)^u}
v_\lambda(Y_v(t)) e^{\lambda (X_u(t) - \xi_{S_u}) - E_\lambda
(t-S_u)}}{\, \tilde{\mathcal{G}}_\infty}.
\end{multline*}
We know from the decomposition \eqref{q_decomp} that the under the
measure $\tilde{\Q}_\lambda$ the subtrees coming off the spine
evolve as if under the measure $P$, and therefore
\begin{gather*}
\tilde{\Q}_\lambda\ccond{\sum_{v \in N_t, v \in (\tau,M)^u}
v_\lambda(Y_v(t)) e^{\lambda (X_u(t) - \xi_{S_u}) - E_\lambda
(t-S_u)}}{\, \tilde{\mathcal{G}}_\infty}\\ =
\tilde{P}\ccond{\sum_{v \in N_t, v \in (\tau,M)^u}
v_\lambda(Y_v(t)) e^{\lambda (X_u(t) - \xi_{S_u}) - E_\lambda
(t-S_u)}}{\, \tilde{\mathcal{G}}_\infty} = v_\lambda(\eta_{S_u}),
\end{gather*}
since the additive expression being evaluated on the subtree is
just a shifted form of the martingale $Z_\lambda$ itself.

This concludes the proof of \eqref{hand}, but before we go move on
to give a similar proof for the general case, for easier reference
through the cumbersome-looking general proof it is worth recalling
that
\[
\zeta_\lambda(t) = e^{\int_0^t R(\eta_s) \, \diffd s}
v_\lambda(\eta_t) e^{\lambda \xi_t - E_\lambda t},
\]
and therefore noting that \eqref{hand} can alternatively be
written as
\[
\tilde{\Q}_\lambda \cond{Z_\lambda(t)}{\tilde{\mathcal{G}}_\infty}
= \sum_{u < N_t} e^{-\int_0^{S_u} R(\eta_s) \, \diffd s}
\zeta_\lambda(S_u) + e^{-\int_0^t R(\eta_s) \, \diffd s}
\zeta_\lambda(t).
\]
Also, in the general model we are supposing that each particle $u$
in the spine will give birth to a total of $A_u$ subtrees that go
off from the spine -- the one remaining other offspring is used to
continue the line of descent that makes up the spine. This
explains the appearance of $A_u$ in the general decomposition.
\begin{thm}[Spine decomposition]\label{thm:spine_decomposition}
We have the following \textbf{spine decomposition} for the
additive branching-particle martingale:
\[
\tilde{\Q}^x\cond{Z(t)}{\tilde{\mathcal{G}}_\infty} = \sum_{u <
\xi_t} A_u \, e^{-\int_0^{S_u} m(\xi_s) R(\xi_s) \, \diffd s}
\zeta(S_u) + e^{-\int_0^t m(\xi_s) R(\xi_s) \, \diffd s} \zeta(t).
\]
\end{thm}
\Proof In each sample tree one and only one of the particles alive
at time $t$ is the spine and therefore:
\begin{align*}
Z(t) &= \sum_{u \in N_t} e^{-\int_0^t m(X_u(s)) R(X_u(s)) \,
\diffd s} \zeta_u(t),\\
&= e^{-\int_0^t m(\xi_s) R(\xi_s) \, \diffd s} \zeta(t) + \sum_{u
\in N_t, u \neq \xi_t} e^{-\int_0^t m(X_u(s)) R(X_u(s)) \, \diffd
s} \zeta_u(t).
\end{align*}
The other individuals $\set{u \in N_t, u \neq \xi_t}$ can be
partitioned into subtrees created from fissions along the spine.
That is, each node $u$ in the spine $\xi_t$ (so $u < \xi_t$) has
given birth at time $S_u$ to one offspring node $uj$ (for some $1
\leq j \leq 1+A_u$) that was chosen to continue the spine whilst
the other $A_u$ individuals go off to make the subtrees
$(\tau,M)^u_j$. Therefore,
\begin{equation}\label{tommy}
Z(t) = e^{-\int_0^t m(\xi_s) R(\xi_s) \, \diffd s} \zeta(t) +
\sum_{u < \xi_t} e^{-\int_0^{S_u} m(\xi_s) R(\xi_s) \, \diffd s}
\sum_{\substack{j=1, \ldots ,1+A_u\\ uj \notin \xi}}
Z_{uj}(S_u;t),
\end{equation}
where for $t \geq S_u$,
\[
Z_{uj}(S_u;t) := \sum_{v \in N_t, v \in (\tau,M)^u_j}
e^{-\int_{S_u}^t m(X_v(s)) R(X_v(s)) \, \diffd s} \zeta_v(t),
\]
is, conditional on $\tilde{\mathcal{G}}_\infty$, a
$\tilde{P}$-martingale on the subtree $(\tau,M)^u_j$, and
therefore
\[
\tilde{P}\cond{Z_{uj}(S_u;t)}{\tilde{\mathcal{G}}_\infty} =
\zeta(S_u).
\]
Thus taking $\tilde{\Q}$-conditional expectations of \eqref{tommy}
gives
\begin{align*}
\tilde{\Q}^x\cond{Z(t)}{\tilde{\mathcal{G}}_\infty} &=
e^{-\int_0^t m(\xi_s) R(\xi_s) \, \diffd s} \zeta(t) + \sum_{u <
\xi_t} e^{-\int_0^{S_u} m(\xi_s) R(\xi_s) \, \diffd s}
\tilde{P}\ccond{\sum_{\substack{j=1, \ldots ,1+A_u\\ uj
\notin \xi}} \!\!\!\!\!\!\!\!Z_{uj}(S_u;t)}{\, \tilde{\mathcal{G}}_\infty},\\
&= e^{-\int_0^t m(\xi_s) R(\xi_s) \, \diffd s} \zeta(t) + \sum_{u
< \xi_t} e^{-\int_0^{S_u} m(\xi_s) R(\xi_s) \, \diffd s} A_u \,
\zeta(S_u),
\end{align*}
which completes the proof. \qed

This representation was first used in the Lyons \emph{et al}
\cite{lyons:simple_biggins_convergence,
lyons_kurtz_pemantle_peres:conceptual_Kesten_Stigum,
lyons_pemantle_peres:conceptual_LlogL} papers and has become the
standard way to investigate the behaviour of $Z$ under the measure
$\tilde{\Q}$. We observe that the two measures $\tilde{P}$ and
$\tilde{\Q}$ for the general model are equal when conditioned on
$\tilde{\mathcal{G}}_\infty$ since this factors out their
differences in the spine diffusion $\xi_t$, the family sizes born
from the spine and the fission times on the spine. Therefore it
follows that same argument as used above applies for $\tilde{P}$
to give:
\begin{cor}
\[
\tilde{P}\cond{Z(t)}{\tilde{\mathcal{G}}_\infty} = \sum_{u <
\xi_t} A_u \, e^{-\int_0^{S_u} m(\xi_s) R(\xi_s) \, \diffd s}
\zeta(S_u) + e^{-\int_0^t m(\xi_s) R(\xi_s) \, \diffd s} \zeta(t).
\]
\end{cor}

\section{Spine results}\label{sec:spine_methods}
Having covered the formal basis for our spine approach, we now
present some new results that follow from our spine formulation:
the Gibbs-Boltzmann weights, conditional expectations, and a
simpler proof of the improved Many-to-One theorem.

\subsection{The \emph{Gibbs-Boltzmann} weights of
$\tilde{\Q}$}\label{sec:g-b_weights}

The Gibbs-Boltzmann weightings in branching processes are
well-known, for example see Chauvin and Rouault
\cite{chauvin_rouault:gibbs_boltzmann} where they consider random
measures on the boundary of the tree, and Harris
\cite{harris:gibbs_boltzmann} which gives convergence results for
Gibbs-Boltzmann random measures. They have previously been
considered via the individual terms of the additive martingale
$Z$, but the following theorem gives a new interpretation of these
weightings in terms of the spine. We recall that
\[
Z(t) = \sum_{u \in N_t} e^{-\int_0^t m(X_u(s)) R(X_u(s)) \, \diffd
s} \zeta_u(t).
\]
\begin{thm}\label{thm_spine_is_particular_node}
Let $u \in \Omega$ be a given and fixed label. Then
\[
\tilde{\Q}\cond{\xi_t = u}{\F_t} = \Ind_{(u\in
N_t)}\frac{e^{-\int_0^t m(X_u(s)) R(X_u(s)) \, \diffd s}
\zeta_u(t)}{Z(t)}.
\]
\end{thm}
\Proof Suppose $u\in\Omega$, and $F\in\F_t$. We aim to show:
\[
\int_F \Ind_{(\xi_t=u)} \,\diffd\tilde{\Q}(\tau,M,\xi) = \int_F
\Ind_{(u\in N_t)}\frac{e^{-\int_0^t m(X_u(s)) R(X_u(s)) \, \diffd
s} \zeta_u(t)}{Z(t)} \,\diffd \tilde{\Q}(\tau,M,\xi).
\]
First of all we know that $\diffd \tilde{\Q}/ \diffd \tilde{P} =
\tilde{\zeta}(t)$ on $\F_t$ and therefore,
\[
\text{LHS} = \int_F \Ind_{(\xi_t=u)} \prod_{v < \xi_t} (1+A_v)
e^{-\int_0^t m(\xi_s)R(\xi_s) \, \diffd s} \zeta(t) \,\,\diffd
\tilde{P}(\tau,M,\xi),
\]
by definition of $\tilde{\zeta}(t)$ at \eqref{eqn:sp_mg_tilde}.
The definition \ref{defn_of_P_tilde} of the measure $\tilde{P}$
requires us to express the integrand with a representation like
\eqref{lyons_f_tilde_representation}:
\begin{gather*}
\Ind_{(\xi_t=u)} \prod_{v < \xi_t} (1+A_v) e^{-\int_0^t
m(\xi_s)R(\xi_s) \, \diffd s} \zeta(t) \\
\begin{aligned}
&= \Ind_{(\xi_t = u)} \sum_{w \in N_t} \prod_{v < w} (1+A_v)
e^{-\int_0^t m(X_w(s))R(X_w(s)) \, \diffd s} \zeta_w(t)
\Ind_{(\xi_t=w)},\\
&= \Ind_{(\xi_t = u)} \Ind_{(u \in N_t)} \prod_{v < u} (1+A_v)
e^{-\int_0^t m(X_u(s))R(X_u(s)) \, \diffd s} \zeta_u(t),
\end{aligned}
\end{gather*}
and therefore
\begin{align*}
\text{LHS} &= \int_F \Ind_{(u \in N_t)} \prod_{v < u} (1+A_v)
e^{-\int_0^t m(X_u(s))R(X_u(s)) \, \diffd s} \zeta_u(t)
\Ind_{(\xi_t = u)} \,\,\diffd
\tilde{P}(\tau,M,\xi),\\
&= \int_F \Ind_{(u \in N_t)} e^{-\int_0^t m(X_u(s))R(X_u(s)) \,
\diffd s} \zeta_u(t) \,\,\diffd P(\tau,M,\xi),
\end{align*}
by definition \ref{defn_of_P_tilde}. We emphasize that now this is
an integral taken with respect to the measure $P$ over the
$\sigma$-algebra $\F_t$, and here we know that $\diffd P/\diffd \Q
= 1/Z(t)$, so:
\[
\text{LHS} = \int_F \Ind_{(u \in N_t)} e^{-\int_0^t
m(X_u(s))R(X_u(s)) \, \diffd s} \zeta_u(t) \,\,\frac{1}{Z(t)}
\diffd \Q(\tau,M,\xi),
\]
and the proof is concluded.\qed

The above result combines with the representation
\eqref{lyons_f_tilde_representation} to show how we take
conditional expectations under the measure $\tilde{\Q}$.
\begin{thm}\label{thm:q_cond_exps}
If $f(t) \in m\tilde{\F}_t$, and $f = \sum_{u\in N_t} f_u(t)
\Ind_{(\xi_t = u)}$, with $f_u(t) \in m\F_t$ then
\begin{equation}
\tilde{\Q}\cond{f(t)}{\F_t} = \sum_{u\in N_t} f_u(t)
\frac{e^{-\int_0^t m(X_u(s)) R(X_u(s)) \, \diffd s}
\zeta_u(t)}{Z(t)}.
\end{equation}
\end{thm}
\Proof It is clear that
\[
\tilde{\Q}\cond{f(t)}{\F_t} = \sum_{u\in N_t} f_u(t)
\tilde{\Q}\cond{\xi_t = u}{\F_t},
\]
and the result follows from Theorem \ref{thm_spine_is_particular_node}.\qed\\
A corollary to this exceptionally useful result also appears to go a
long way towards obtaining the Kesten-Stigum result in more general
models:
\begin{cor}
If $g(\cdot)$ is a Borel function on $J$ then
\begin{equation}\label{kslink}
\sum_{u\in N_t} g(X_u(t)) \, e^{-\int_0^t m(X_u(s)) R(X_u(s)) \,
\diffd s} \zeta_u(t) = \tilde{\Q}\cond{g(\xi_t)}{\F_t} \times
Z(t).
\end{equation}
\end{cor}
\Proof We can write $g(\xi_t) = \sum_{u\in N_t}
g(X_u(t)) \Ind{(\xi_t = u)}$, and now the result follows from the
above corollary.\qed

\noindent The classical Kesten-Stigum theorems of
\cite{kesten_stigum_1,kesten_stigum_2,kesten_stigum_3} for
multi-dimensional Galton-Watson processes give conditions under
which an operation like the left-hand side of \eqref{kslink}
converges as $t \to \infty$, and it is found that when it exists the
limit is a multiple of the martingale limit $Z(\infty)$. Also see
Lyons \emph{et al}
\cite{lyons_kurtz_pemantle_peres:conceptual_Kesten_Stigum} for a
recent proof of this based on other spine techniques. Our spine
formulation apparently gives a previously unknown but simple meaning
to this operation in terms of a conditional expectation and, as we
hope to pursue in further work, in many cases we would intuitively
expect that
$\tilde{\Q}\cond{g(\xi_t)}{\F_t}/\tilde{\Q}(g(\xi_t))\rightarrow 1$
a.s., leading to alternative spine proofs of both Kesten-Stigum like
theorems and Watanabe's theorem in the case of BBM.

\subsection{The \emph{Full} Many-to-One
Theorem}\label{sec:full_many-to-one} An very useful tool in the
study of branching processes is the \emph{Many-to-One} result that
enables expectations of sums over particles in the branching process
to be calculated in terms of an expectation of a single particle. In
the context of the finite-type branching diffusion of section
\ref{sec:finite-type_model}, the Many-to-One theorem would be stated
as follows:
\begin{thm}
For any measurable function $f: J \to \Real$ we have
\[
P^{x,y}\bseq{\sum_{u \in N_t} f(X_u(t),Y_u(t))} = \Prob^{x,y}
\bseq{e^{\int_0^t R(\eta_s)\,\diffd s} f(\xi_t,\eta_t)}.
\]
\end{thm}
Intuitively it is clear that the up-weighting term $e^{\int_0^t
R(\eta_s)\,\diffd s}$ incorporates the notion of the population
growing at an exponential rate, whilst the idea of
$f(\xi_t,\eta_t)$ being the `typical' behaviour of
$f(X_u(t),Y_u(t))$ is also reasonable.

Existing results tend to apply only to functions of the above form
that depend only on \emph{the time-$t$ location} of the spine and
existing proofs do not lend themselves to covering functions that
depend on the entire \emph{path history} of the spine up to time
$t$.

With the spine approach we have the benefit of being able to give
a much less complicated proof of the stronger version that covers
the most general path-dependent functions.
\begin{thm}[Many-to-One]\label{thm:full_many-to-one}
If $f(t) \in m\tilde{\mathcal{F}}_t$ has the representation
\[
f(t) = \sum_{u \in N_t} f_u(t) \Ind_{(\xi_t = u)},
\]
where $f_u(t) \in m\F_t$, then
\begin{equation}\label{M21general}
P\bseq{\sum_{u\in N_t} f_u(t) e^{-\int_0^t m(X_u(s)) R(X_u(s)) \,
\diffd s} \zeta_u(t)}
= \tilde{P}\seq{f(t)\, \tilde{\zeta}(t)}
=  {\zeta}(0)\,\tilde\Q\seq{f(t)}.
\end{equation}
In particular, if $g(t) \in m\mathcal{G}_t$ with %has the representation
$
g(t) = \sum_{u \in N_t} g_u(t) \Ind_{(\xi_t = u)}
$
where $g_u(t) \in m\F_t$, then
\begin{equation}\label{M21specialcase}
P\bseq{\sum_{u\in N_t} g_u(t) }
= \Prob\bseq{e^{\int_0^t m(\xi_s)R(\xi_s) \diffd s} g(t)}
= \hat\Prob\bseq{\frac{g(t)\,\zeta(0)}{e^{-\int_0^t m(\xi_s)R(\xi_s) \diffd s}\,\zeta(t)}}.
\end{equation}
\end{thm}
\Proof Let $f(t)$ be an $\tilde{\F}_t$-measurable function with the
given representation. We can use the tower property together with
Theorem \ref{thm:q_cond_exps} to obtain
\begin{align*}
\tilde{\Q}\seq{f(t)}
&=\tilde{\Q}\bseq{\tilde{\Q}\cond{f(t)}{\F_t}}
= \Q\bseq{\tilde{\Q}\cond{f(t)}{\F_t}} \\
&= \Q\bseq{\frac{1}{Z(t)} \sum_{u\in N_t} f_u(t) e^{-\int_0^t
m(X_u(s)) R(X_u(s)) \, \diffd s} \zeta_u(t)}.
\end{align*}
We emphasize that this is a $\Q$ expectation of a
$\F_t$-measurable expression. From Theorem \ref{mark_riley},
\[
\frac{\diffd \Q}{\diffd P}\bigg\vert_{\F_t} = \frac{Z(t)}{Z(0)},
\]
and therefore we have
\[
\tilde{\Q}\seq{f(t)} = P\bseq{ Z(0)^{-1}\,\sum_{u\in N_t} f_u(t)
e^{-\int_0^t m(X_u(s)) R(X_u(s)) \, \diffd s} \zeta_u(t)}.
\]

On the other hand, since $f(t)$ is $\tilde\F_t$-measurable and
\[
\frac{\diffd \tilde{\Q}}{\diffd
\tilde{P}}\bigg\vert_{\tilde{\F}_t} = \frac{\tilde{\zeta}(t)}{\tilde{\zeta}(0)},
\]
we have
\[
\tilde{\Q}\seq{f(t)} =  \tilde{P}\seq{f(t) \times
\tilde{\zeta}(t)\,\tilde{\zeta}(0)^{-1}}.
\]
Trivially noting $Z(0)=\zeta(0)=\tilde\zeta(0)$ as there is only one
initial ancestor, we can combine these expressions to obtain
\eqref{M21general}. For the second part, given $g(t) \in
m\mathcal{G}_t$, we can define
\[
f(t) := e^{\int_0^t m(\xi_s) R(\xi_s) \, \diffd s} g(t) \times
\zeta(t)^{-1},
\]
which is clearly $\mathcal{G}_t$-measurable and satisfies $f(t) =
\sum_{u\in N_t}  f_u(t) \Ind_{(\xi_t = u)}$ with
\[
f_u(t) = g_u(t) e^{\int_0^t m(X_u(s)) R(X_u(s)) \, \diffd s}
\zeta_u(t)^{-1} \, \in m\F_t.
\]
When we use this $f(t)$ in equation \eqref{M21general}
and recall Lemma \ref{z_and_zeta_from_zeta_tilde}, that
${\Prob} := \tilde{P}\vert_{\mathcal{G}_\infty}$ from Definition \ref{dave}
and that $\hat{\Prob} := \tilde{\Q}\vert_{\mathcal{G}_\infty}$ from Definition \ref{dave2},
we arrive at the particular case given at \eqref{M21specialcase} in the theorem. \qed

In the further special case in which $g = g(\xi_t)$ for some Borel-measurable
function $g(\cdot)$, the trivial representation
\[
g(\xi_t) = \sum_{u \in N_t} g\seq{X_u(t)}
\]
leads immediately to the weaker version of the Many-to-One result
that was utilised and proven, for example, in
Harris and Williams \cite{art:largedevs1} and Champneys \emph{et al} \cite{art:aap}
using resolvents and the Feynman-Kac formula, expressed in terms of our more
general branching Markov process $\mathbb{X}_t$:
\begin{cor}\label{thm:many-to-one_in_spine_chapter}
If $g(\cdot): J \to \Real$ is $\mathcal{B}$-measurable then
\[
P\bseq{\sum_{u\in N_t} g(X_u(t))} = \Prob\bseq{e^{\int_0^t
R(\xi_s) \diffd s} g(\xi_t)}.
\]
\end{cor}

% ----------------------------------------------------------------
%
%
%
\def\polhk#1{\setbox0=\hbox{#1}{\ooalign{\hidewidth
  \lower1.5ex\hbox{`}\hidewidth\crcr\unhbox0}}} \def\cprime{$'$}
\providecommand{\bysame}{\leavevmode\hbox
to3em{\hrulefill}\thinspace}
\providecommand{\MR}{\relax\ifhmode\unskip\space\fi MR }
\providecommand{\MRhref}[2]{%
  \href{http://www.ams.org/mathscinet-getitem?mr=#1}{#2}
} \providecommand{\href}[2]{#2}

\end{document}